\documentclass[12pt,reqno]{amsart}
\usepackage{amscd,amsmath,amsthm,amssymb}
\usepackage{color}
\usepackage{pstricks}
\usepackage{stmaryrd}
\usepackage{url}

\usepackage{tikz-cd}
\newcommand{\arrowIn}{
	\tikz \draw[-stealth] (-1pt,0) -- (1pt,0);
}
\newcommand{\arrowOut}{
	\tikz \draw[-stealth] (1pt,0) -- (-1pt,0);
}
\usepackage{latexsym}
\usepackage{amsfonts,amsmath,mathtools}
\usepackage{graphics}
\usepackage{float}
\usepackage{enumitem}

\usepackage{booktabs} 
\usepackage{colortbl}
\usepackage{lipsum}

\newpsstyle{fatline}{linewidth=1.5pt}
\newpsstyle{fyp}{fillstyle=solid,fillcolor=verylight}
\definecolor{verylight}{gray}{0.97}
\definecolor{light}{gray}{0.9}
\definecolor{medium}{gray}{0.85}
\definecolor{dark}{gray}{0.6}

\def\NZQ{\mathbb}               

\def\ZZ{{\NZQ Z}}

\def\G{{\mathcal G}}

\def\pd{\textup{proj}\phantom{.}\!\textup{dim}}

\def\opn#1#2{\def#1{\operatorname{#2}}} 
\opn\chara{char} \opn\length{\ell} \opn\pd{pd} \opn\rk{rk}
\opn\projdim{proj\,dim} \opn\injdim{inj\,dim} \opn\rank{rank}
\opn\depth{depth} \opn\grade{grade} \opn\height{height}
\opn\embdim{emb\,dim} \opn\codim{codim}

\opn\Tr{Tr} \opn\bigrank{big\,rank}
\opn\superheight{superheight}\opn\lcm{lcm}
\opn\trdeg{tr\,deg}
	\opn\reg{reg} \opn\lreg{lreg} \opn\ini{in} \opn\lpd{lpd}
	\opn\size{size} \opn\sdepth{sdepth}
	\opn\link{link}\opn\fdepth{fdepth}\opn\lex{lex}
	\opn\tr{tr}
	\opn\type{type}
	\opn\gap{gap}
	\opn\diam{diam}
	\opn\Mod{Mod}
	\opn\div{div} \opn\Div{Div} \opn\cl{cl} \opn\Cl{Cl}
	\opn\Spec{Spec} \opn\Supp{Supp} \opn\supp{supp} \opn\Sing{Sing}
	\opn\Ass{Ass} \opn\Min{Min}\opn\Mon{Mon}
	\opn\Ann{Ann} \opn\Rad{Rad} \opn\Soc{Soc}
	\opn\Im{Im} \opn\Ker{Ker} \opn\Coker{Coker} \opn\Am{Am}
	\opn\Hom{Hom} \opn\Tor{Tor} \opn\Ext{Ext} \opn\End{End}
	\opn\Aut{Aut} \opn\id{id}
	
	\opn\nat{nat}
	\opn\pff{pf}
	\opn\Pf{Pf} \opn\GL{GL} \opn\SL{SL} \opn\mod{mod} \opn\ord{ord}
	\opn\Gin{Gin} \opn\Hilb{Hilb}\opn\sort{sort}
	\opn\PF{PF}\opn\Ap{Ap}
	\opn\dist{dist}
	\opn\aff{aff}
	\opn\relint{relint} \opn\st{st}
	\opn\lk{lk} \opn\cn{cn} \opn\core{core} \opn\vol{vol}  \opn\inp{inp} \opn\nilpot{nilpot}
	\opn\link{link} \opn\star{star}\opn\lex{lex}\opn\set{set}
	\opn\width{wd}
	\opn\Fr{F}
	\opn\QF{QF}
	\opn\G{G}
	\opn\type{type}\opn\res{res}
	\opn\conv{conv}
	\opn\sr{sr}
	\opn\gr{gr}
	
	\def\pot#1#2{#1[\kern-0.28ex[#2]\kern-0.28ex]}

	\opn\dirlim{\underrightarrow{\lim}}
	\opn\inivlim{\underleftarrow{\lim}}
	\def\Implies{\ifmmode\Longrightarrow \else
		\unskip${}\Longrightarrow{}$\ignorespaces\fi}
	\def\implies{\ifmmode\Rightarrow \else
		\unskip${}\Rightarrow{}$\ignorespaces\fi}
	\def\iff{\ifmmode\Longleftrightarrow \else
		\unskip${}\Longleftrightarrow{}$\ignorespaces\fi}

	\let\:=\colon
	\newtheorem{Theorem}{Theorem}[section]
	
	\newtheorem{Corollary}[Theorem]{Corollary}

	\let\epsilon\varepsilon
	\let\kappa=\varkappa
	\def\qed{\ifhmode\textqed\fi
		\ifmmode\ifinner\hfill\quad\qedsymbol\else\dispqed\fi\fi}
	\def\textqed{\unskip\nobreak\penalty50
		\hskip2em\hbox{}\nobreak\hfill\qedsymbol
		\parfillskip=0pt \finalhyphendemerits=0}
	\def\dispqed{\rlap{\qquad\qedsymbol}}
	
	\opn\dis{dis}
	\def\pnt{{\raise0.5mm\hbox{\large\bf.}}}
	
	\opn\Lex{Lex}
	\opn\Max{Max}
	\opn\Shad{Shad}
	\opn\astab{astab}

	\opn\v{v}

\begin{document}
	
\title{Matching Powers: Macaulay2 Package}	
\author{Antonino Ficarra}

\address{Antonino Ficarra, Department of mathematics and computer sciences, physics and earth sciences, University of Messina, Viale Ferdinando Stagno d'Alcontres 31, 98166 Messina, Italy}
\email{antficarra@unime.it}

\thanks{
}



\maketitle

\begin{abstract}
	We introduce the \textit{Macaulay2} package \texttt{MatchingPowers}. It allows to compute and manipulate the matching powers of a monomial ideal. The basic theory of matching powers is explained and the main features of the package are presented.
\end{abstract}

\section{Basic Properties of Matching Powers}

Let $S=K[x_1,\dots,x_n]$ be the standard graded polynomial ring with coefficients over a field $K$. Let $I$ be a monomial ideal of $S$. The $k$th \textit{matching power} of $I$, denoted by $I^{[k]}$ is the monomial ideal generated by the products $u_1\cdots u_k$ such that $u_1,\dots,u_k$ is a monomial regular sequence of length $k$ contained in $I$. The biggest $k$ such that $I^{[k]}\ne0$ is called the \textit{monomial grade} of $I$ and it is denoted by $\nu(I)$. The concept of matching powers was first introduced in the squarefree case in \cite{EHHM2022b} and later extended to all monomial ideals in \cite{EF}.

One of the main motivations for considering matching powers comes from graph theory. Let $G$ be a finite simple graph on vertex set $[n]=\{1,\dots,n\}$. The \textit{edge ideal} of $G$ is the squarefree monomial ideal $I(G)$ generated by the monomials $x_e=x_ix_j$ such that $e=\{i,j\}$ is an edge of $G$. Let $e_1,\dots,e_k$ be $k$ edges of $G$, then $x_{e_1}\cdots x_{e_k}$ is a generator of $I(G)^{[k]}$ if and only if $\{e_1,\dots,e_k\}$ is a $k$-matching of $G$. This observation justifies the choice to name the powers $I^{[k]}$ the matching powers of $I$.

Similarly, for a non-squarefree example, we can consider edge ideals of \textit{weighted oriented graphs}. Recall that a (\textit{vertex})-\textit{weighted oriented graph} $\mathcal{D}=(V(\mathcal{D}),E(\mathcal{D}),w)$ consists of an underlying graph $G$, with $V(\mathcal{D})=V(G)=[n]$, on which each edge is given an orientation and it is equipped with a \textit{weight function} $w:V(G)\rightarrow\mathbb{Z}_{\ge1}$. The \textit{weight} of a vertex $i\in V(G)$, denoted by $w_i$, is the value $w(i)$ of the weight function at $i$. The directed edges of $\mathcal{D}$ are denoted by pairs $(i,j)\in E(\mathcal{D})$ to reflect the orientation, hence $(i,j)$ represents an edge directed from $i$ to $j$. The \textit{edge ideal} of $\mathcal{D}$ is defined as the ideal
$$
I(\mathcal{D})\ =\ (x_ix_j^{w_j}\ :\ (i,j)\in E(\mathcal{D})).
$$ 

As before, given $k$ generators $u_\ell=x_{i_\ell}x_{j_\ell}^{w_{j_\ell}}$ of $I(\mathcal{D})$, with $e_\ell=(i_\ell,j_\ell)\in E(\mathcal{D})$, for $\ell=1,\dots,k$, the product $u_1\cdots u_k$ is a generator of $I(\mathcal{D})^{[k]}$ if and only if $\{e_1,\dots,e_k\}$ is a $k$-matching of the underlying simple graph $G$ of $\mathcal{D}$.

In the remaining part of this section, we summarize some of the main results proved in \cite[Section 1]{EF}.

Let $G(I)$ be the unique minimal set of monomial generators of $I$. The \textit{initial degree} of $I$ is defined as $\textup{indeg}(I)=\min\{\deg(u):u\in G(I)\}$. For a monomial $u\in S$, the \textit{$x_i$-degree} of $u$ is defined as $\deg_{x_i}(u)=\max\{j:x_i^j\ \textup{divides}\ u\}$. Whereas, the \textit{support} of $u$ is the set $\supp(u)=\{i:\deg_{x_i}(u)>0\}$. Following \cite{F2}, we define the \textit{bounding multidegree} of $I$ as
$$
{\bf deg}(I)\ =\ (\deg_{x_1}(I),\dots,\deg_{x_n}(I)),
$$
with
$$
\deg_{x_i}(I)\ =\ \max_{u\in G(I)}\deg_{x_i}(u),\ \ \textup{for all}\ \ \ 1\le i\le n.
$$

Furthermore, if ${\bf a}=(a_1,\dots,a_n)\in\ZZ_{\ge0}^n$, we set $|{\bf a}|=a_1+\dots+a_n$.

For a monomial $u=x_1^{a_1}\cdots x_n^{a_n}\in S$, we set $u^\wp=\prod_{i=1}^n(\prod_{j=1}^{a_i}x_{i,j})$ in the polynomial ring $K[x_{i,j}:1\le i\le n,1\le j\le a_i]$. The \textit{polarization} of $I$ is defined to be the squarefree ideal $I^\wp$ of $S^\wp=K[x_{i,j}:1\le i\le n,1\le j\le\deg_{x_i}(I)]$ whose minimal generating set is $G(I^\wp)=\{u^\wp:u\in G(I)\}$.

\begin{Theorem}\label{Thm:EF}
	Let $I\subset S$ be a monomial ideal. The following facts hold.
	\begin{enumerate}
		\item[\textup{(i)}] $I^{[k]}=(u_1\cdots u_k\ :\ u_i\in G(I),\ \supp(u_i)\cap\supp(u_j)=\emptyset,\ 1\le i<j\le k)$.
		\item[\textup{(ii)}] $I^{[k]}\ne0$ if and only if $1\le k\le\nu(I)$.
		\item[\textup{(iii)}] $(I^{[k]})^\wp=(I^\wp)^{[k]}$ for all $1\le k\le\nu(I)$.
		\item[\textup{(iv)}] $\depth(S/I^{[k]})=\depth(S^\wp/(I^{\wp})^{[k]})-|{\bf deg}(I)|+n$, for all $1\le k\le\nu(I)$.
		\item[\textup{(v)}] $\depth(S/I^{[k]})\ge\textup{indeg}(I^{[k]})-1+(n-|{\bf deg}(I)|)$.
	\end{enumerate}
\end{Theorem}

Based on inequality (v), we define the \textit{normalized depth function} of $I$ as
$$
g_I(k)\ =\ \depth(S/I^{[k]})+|{\bf deg}(I)|-n-(\textup{indeg}(I^{[k]})-1),
$$
for all $1\le k\le\nu(I)$. In the squarefree case, this is the function introduced in \cite{EHHM2022b}.

As a consequence of Theorem \ref{Thm:EF} it was proved that \cite{EF}
\begin{Corollary}\label{Cor:EF}
	We have $\nu(I)=\nu(I^\wp)$ and $g_I(k)=g_{I^\wp}(k)$ for all $k$.
\end{Corollary}

It was conjectured in \cite{EHHM2022b} that $g_I$ would be non-increasing for all squarefree monomial ideals $I$. Because of Corollary \ref{Cor:EF}, it the above conjecture was true, then $g_I$ would be non-increasing for all monomial ideals $I$ \cite[Proposition 1.8]{EF}. However, very recently, Fakhari \cite{SASF2023b} constructed a family of monomial ideals $I$ such that the difference $g_{I}(2)-g_I(1)$ may be arbitrary large, thus disproving the above conjecture. 

It is an open question whether for any sequence of non-negative integers $a_1,\dots,a_g$, there exists a monomial ideal $I$ such that $\nu(I)=g$ and $g_I(k)=a_k$ for all $k=1,\dots,g$.\medskip

Currently, the study of matching powers of monomial ideals and of various kind of edge ideals \cite{MPV,PS13,PRT} is a very active area of research \cite{BHZN18, CFL, EF, EH2021, EHHM2022a, EHHM2022b, FHH23, SASF2022, SASF2023, SASF2023b}.\smallskip

In the next two sections, we illustrate and explain how to use the \textit{Macaulay2} \cite{GDS} package \texttt{MatchingPowers}. In Section \ref{sec2-FHSPack}, we describe the main functions of the package. In Section \ref{sec3-FHSPack}, an example is presented, that demonstrates how to use the package.

\section{The Package}\label{sec2-FHSPack}

The purpose of the package \texttt{MatchingPowers} is to provide the tools to manipulate and calculate the matching powers of a monomial ideal $I$.

The \textit{multidegree} of a monomial $u\in S$ is defined to be the vector
$$
{\bf deg}(u)\ =\ (\deg_{x_1}(u),\dots,\deg_{x_n}(u)).
$$
Whereas, for ${\bf a}=(a_1,\dots,a_n)\in\ZZ_{\ge0}^n$, we set ${\bf x^a}=\prod_ix_i^{a_i}$. In particular, $u={\bf x}^{{\bf \deg}(u)}$.

Given two monomial ideals $I,J\subset S$, we define the \textit{matching product} of $I$ and $J$ to be the monomial ideal
$$
I*J\ =\ (uv\ :\ u\in I,v\in J,\ u,v\ \textup{is a regular sequence}).
$$
It is not difficult to see that
$$
I*J\ =\ (uv\ :\ u\in G(I),\ v\in G(J),\ \supp(u)\cap\supp(v)=\emptyset).
$$
Moreover, $I^{[2]}=I*I$ and by induction $I^{[k]}=I^{[k-1]}*I$ for all $k\ge2$.

Finally, we recall that a graded ideal $I\subset S$ is \textit{linearly related} if $I$ is generated in single degree, say $d$, and the first syzygy module of $I$ is generated $d+1$.\smallskip

The next table collects the functions available in the package. By $I$ and $J$ we denote two monomial ideals of $S$, ${\bf a}\in\ZZ_{\ge0}^n$ an integral vector, ${\bf x^a}$ a monomial.
\small\begin{table}[H]
	\centering
	\begin{tabular}{ll}
		\rowcolor{black!20}\bottomrule[1.05pt]
		Functions&Description\\
		\toprule[1.05pt]
		\texttt{toMonomial$(S,{\bf a})$}&Computes the monomial ${\bf x^a}$ if ${\bf a}\in\ZZ^n_{\ge0}$\\
		\texttt{toMultidegree}$({\bf x^a})$&Computes the multidegree ${\bf a}$ of ${\bf x^a}$\\
		\texttt{boundingMultidegree}$(I)$&Computes $\textbf{deg}(I)$\\
		\texttt{matchingProduct}$(I,J)$&Computes the matching product of $I$ and $J$\\
		\texttt{matchingPower}$(I,k)$&Computes the $k$th matching power of $I$\\
		\texttt{monomialGrade}$(I)$&Computes the monomial grade $\nu(I)$\\
		\texttt{gFunction}$(I)$&Computes the normalized depth function $g_I$\\
		\texttt{isLinearlyRelated}$(I)$&Checks if $I$ is linearly related\\
		\bottomrule[1.05pt]
	\end{tabular}\medskip
	\caption{List of the functions of \texttt{MatchingPowers}.}
\end{table}
\normalsize

\section{An Example}\label{sec3-FHSPack}

Consider the weighted oriented graph $\mathcal{D}$ depicted below.\bigskip

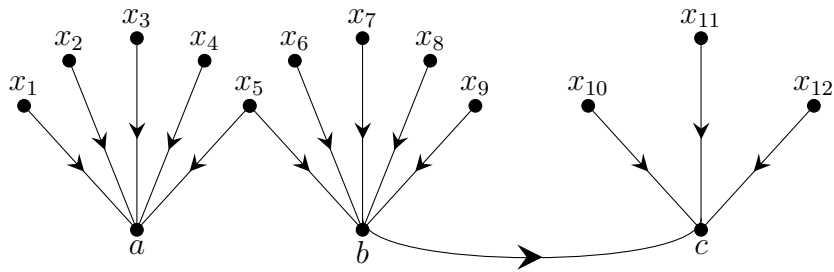
\begin{center}
	\begin{tikzpicture}[scale=1.5]
		\tikzcdset{arrow style=tikz}
		\filldraw (0,0) circle (1.6pt) node[below]{$a$};
		\filldraw (0.6,1.5) circle (1.6pt) node[above]{$x_4$};
		\filldraw (0,1.7) circle (1.6pt) node[above]{$x_3$};
		\filldraw (-0.6,1.5) circle (1.6pt) node[above]{$x_2$};
		\filldraw (-1,1.1) circle (1.6pt) node[above]{$x_1$};
		\filldraw (1,1.1) circle (1.6pt) node[above]{$x_5$};
		\draw (0,1.7) -- (0,0) node[sloped,pos=0.5,xscale=2,yscale=2]{\arrowIn};
		\draw (-0.6,1.5) -- (0,0) node[sloped,pos=0.5,xscale=2,yscale=2]{\arrowIn};
		\draw (0,0) -- (0.6,1.5) node[sloped,pos=0.5,xscale=2,yscale=2]{\arrowOut};
		\draw (-1,1.1) -- (0,0) node[sloped,pos=0.5,xscale=2,yscale=2]{\arrowIn};
		\draw (1,1.1) -- (0,0) node[sloped,pos=0.5,xscale=2,yscale=2]{\arrowOut};
		\filldraw (2,0) circle (1.6pt) node[below]{$b$};
		\filldraw (2,1.7) circle (1.6pt) node[above]{$x_7$};
		\filldraw (1.4,1.5) circle (1.6pt) node[above]{$x_6$};
		\filldraw (1,1.1) circle (1.6pt);
		\filldraw (2.6,1.5) circle (1.6pt) node[above]{$x_8$};
		\filldraw (3,1.1) circle (1.6pt) node[above]{$x_9$};
		\draw (2,1.7) -- (2,0) node[sloped,pos=0.5,xscale=2,yscale=2]{\arrowIn};
		\draw (1.4,1.5) -- (2,0) node[sloped,pos=0.5,xscale=2,yscale=2]{\arrowIn};
		\draw (1,1.1) -- (2,0) node[sloped,pos=0.5,xscale=2,yscale=2]{\arrowIn};
		\draw (2.6,1.5) -- (2,0) node[sloped,pos=0.5,xscale=2,yscale=2]{\arrowOut};
		\draw (3,1.1) -- (2,0) node[sloped,pos=0.5,xscale=2,yscale=2]{\arrowOut};
        \node (a) at (2,0.2) {};
        \node (b) at (5,0.2) {};
        \draw[-] (a)  to [out=-90,in=-90, looseness=0.4] (b) node[xshift=-63,yshift=-19,xscale=3,yscale=3]{\arrowIn};
        \filldraw (5,0) circle (1.6pt) node[below]{$c$};
        \filldraw (5,1.7) circle (1.6pt) node[above]{$x_{11}$};
        \filldraw (4,1.1) circle (1.6pt) node[above]{$x_{10}$};
        \filldraw (6,1.1) circle (1.6pt) node[above]{$x_{12}$};
        \draw (5,1.7) -- (5,0) node[sloped,pos=0.5,xscale=2,yscale=2]{\arrowIn};
        \draw (4,1.1) -- (5,0) node[sloped,pos=0.5,xscale=2,yscale=2]{\arrowIn};
        \draw (6,1.1) -- (5,0) node[sloped,pos=0.5,xscale=2,yscale=2]{\arrowOut};
	\end{tikzpicture}
\end{center}

We assign the weights as follows. For $i=1,\dots,12$, we let $w(x_i)=1$, whereas we set $w(a)=w(b)=w(c)=3$. Thus the edge ideal of $\mathcal{D}$ is
$$
I(\mathcal{D})=(x_1,x_2,x_3,x_4,x_5)(a^3)+(x_5,x_6,x_7,x_8,x_9)(b^3)+(b,x_{10},x_{11},x_{12})(c^3).
$$

Then, $\nu(I(\mathcal{D}))=3$, and by using the package we checked that $I(\mathcal{D})^{[3]}$ is linearly related. Indeed, this also follows from \cite[Theorem 4.2]{EF}.\bigskip

\verb| i1:  S = QQ[x_1..x_12,a..c];|\smallskip

\verb| i2:  I = ideal(x_1,x_2,x_3,x_4,x_5)*ideal(a^3)|\smallskip

\verb|      + ideal(x_5,x_6,x_7,x_8,x_9)*ideal(b^3)|\smallskip

\verb|      + ideal(b,x_10,x_11,x_12)*ideal(c^3);|\smallskip

\texttt{\phantom{i}i3: loadPackage "MatchingPowers"}\smallskip

\texttt{\phantom{i}i4: matchingPower(I,1) == I}\smallskip

\texttt{\phantom{i}o4: true}\smallskip

\texttt{\phantom{i}i5: isLinearlyRelated I}\smallskip

\texttt{\phantom{i}o5: false}\smallskip

\texttt{\phantom{i}i6: J = matchingPower(I,3);}\smallskip

\texttt{\phantom{i}i7: isLinearlyRelated J}\smallskip

\texttt{\phantom{i}o7: true}\smallskip

\texttt{\phantom{i}i8: monomialGrade I}\smallskip

\texttt{\phantom{i}o8: 3}\bigskip

\noindent
A good companion to the \texttt{MatchingPowers} package is the \texttt{HomologicalShiftIdeals} package \cite{FPack1}. Using this latter package, one can check that $I(\mathcal{D})^{[3]}$ is a polymatroidal ideal, which again follows from \cite[Theorem 4.2]{EF}.

\end{document}